\documentclass[12pt,a4paper]{amsart}

\usepackage{fullwidth}
\usepackage[top=1.5in, bottom=1.5in, left=1in, right=1in]{geometry}
\usepackage{amsgen}
\usepackage{amsmath}
\usepackage{amstext}
\usepackage{amsbsy}
\usepackage{amsopn}
\usepackage{amsfonts}
\usepackage{amssymb}
\usepackage{comment}
\usepackage{eepic}
\usepackage{graphicx}
\usepackage{epsf}
\usepackage{pstricks}

\title{A tribute to Mario Petrich}

\author[P.\ V.\ Silva]{Pedro\ V.\ Silva}
\address{Department of Mathematics, University of Porto, Porto, Portugal}
\email{pvsilva@fc.up.pt}

\author[M.\ B.\ Szendrei]{M\'{a}ria\ B.\ Szendrei}
\address{Bolyai Institute, University of Szeged, Szeged, Hungary}
\email{m.szendrei@math.u-szeged.hu}

\begin{document}

$\ $

\vskip -1.5in

\noindent
{\tiny%\small
{\hfill This preprint has not undergone peer review or any post-submission improvements or}

{\hfill corrections. The article is to be published in Semigroup Forum.}}%, and is available online at}

%{\hfill {\tt https://doi.org/10.1007/s00233-021-10235-5}.}}

\vskip 1.2in

\maketitle

Mario Petrich has had a great influence on semigroup theory and on several generations of semigroup theorists. 
His extensive and deep contributions to the theory of regular, in particular, completely regular and inverse semigroups, and to a great variety of other subjects, earned him a special place in the history of our relatively young area.
He is appreciated and admired by many of the people who have got the chance to know him and have benefitted from his insight and ideas.

Mario Petrich was born in 1932 in Split, at that time in Yugoslavia. 
He spent part of his childhood and adolescence there under Italian rule, during World War II. 
Then Split was recovered by Yugoslavia, and he lived there until he was 20 years old, completing three years of studies in mathematics at the University of Zagreb. 
Then he left Yugoslavia, first went to West Germany and then to the USA. 
There he first worked as a messenger boy in New York at \$1 per hour for 7 months, and then continued his studies. 
He completed his M.Sc.\ at Midwest University, and took his Ph.D.\ at the University of Washington, Seattle in 1961. 
After working at various institutions for shorter terms, partly in applied mathematics, he took a position at Pennsylvania State University in 1964. 
In 1976 he left the United States for personal reasons. 
He became a `citizen of the world', accepting limited term appointments. 
Over the years, he was hosted, as a visiting professor or such, by many universities in various countries where he engaged in many successful collaborations: 
in Austria, Canada, France, Germany, Italy, Portugal, Spain, the UK. 
Finally, he settled in a village on an island on the Dalmatian coast, Croatia, continuing productive research in semigroup theory. 
His range of interests extends well beyond mathematics, going from the intellectual-cultural sphere to physical activities (a fierce swimmer and hiker).

Mario has fundamental results on various classes of semigroups. 
So far, he has published more than 260 articles, mainly on semigroups, but he has papers, for example, also on rings, lattices and projective geometry.
Among his semigroup papers, over 80 are connected with completely regular semigroups, about 30 with inverse semigroups, and about one half of the rest with other types of regular semigroups.
So it is hard to do justice to his publications, it is only possible to touch on some highlights and themes.

From the first period of research of Mario Petrich, his work on semilattice decompositions \cite{Pet1964} should be mentioned. 
He introduced a new approach to semilattice decompositions using completely prime ideals (notice that in this paper he uses, following the terminology in \cite{CliPre1961-7}, the old (by now outdated) name `prime ideal' for them). 
Using this presentation, he also characterizes several classes of semigroups such that all classes of their greatest semilattice decomposition have a given property.  
His work in this direction has been rounded up in the book \cite[Chapter II]{Pet1973m}. 
In the opposite direction we have the \emph{semilattice composition problem}: 
how a semigroup $S$ can be reconstructed from a semilattice decomposition 
$S= \bigcup_{i\in Y} S_i$ if certain mappings and conditions on them are given which help to define multiplication between elements from different components. 
This problem was solved for completely regular semigroups (that is, for unions of groups) by Lallement \cite{Lall1967} and in the general case by Petrich, published as well in the book \cite{Pet1973m}, see also \cite[Theorem III.2.8]{Gril1995}. 

A part of Mario's research has been greatly motivated by classical results on rings.
An outstanding theory developed by him on dense (ideal) extensions of completely $0$-simple semigroups (\cite{Pet1973.1},\cite{Pet1973.2}) generalizes a number of results achieved formerly on semigroups and rings which did not appear to have much in common.
The Jacobson theorem on the structure of primitive rings with nonzero socle 
\cite[Section IV.9]{Jac1956} is among these results, and its proof has served as a model for approaching the semigroup case. 
Up to isomorphism, this theorem characterizes the rings mentioned as the rings of linear transformations associated in a way with pairs of dual vector spaces over division rings.
A semigroup analogue of this theory was developed by Hoehnke \cite{Hoe1966} to find representations of certain primitive semigroups. 
Given a group with zero $G^0$, Mario defines the notions of a {\em weakly dual pair of exact $G^0$-acts}.
For every weakly dual pair $(R,L)$ he considers the semigroup of adjoints $\mathcal{E}(R,L)$ of $(R,L)$, and a dense completely $0$-simple ideal $\mathcal{C}(R,L)$ in it. 
One of the main results of Mario's two-piece paper is that a semigroup is a dense extension of a completely $0$-simple semigroup if and only if it is isomorphic to a subsemigroup of a semigroup $\mathcal{E}(R,L)$ for some weakly dual pair $(R,L)$ which contains $\mathcal{C}(R,L)$.
Moreover, a representation of the same abstract semigroups is provided by pairs of partial transformations.   
If the completely $0$-simple semigroup is (left) reductive then the weakly dual pair of exact $G^0$-acts can be replaced by a single exact $G^0$-act, and the pair of partial transformations by a single partial transformation.
The isomorphisms between semigroups of adjoints are described by means of semi-isomorphisms between the 
underlying weakly dual pairs, where a semi-isomorphism is an analogue of a semilinear isomorphism with an adjoint.    
The main results of the paper can be applied to derive a number of structure and representation theorems on semigroups, see for example  Gluskin \cite{Glu1959a} and Ljapin \cite{Lja1955},\cite{Lja1964}.
Moreover, they allow Mario to add several new necessary and sufficient conditions, partly motivated by Gluskin \cite{Glu1959b} and Lawson \cite{LawL1969}, to the Jacobson theorem, one of them requiring that the multiplicative semigroup of the ring is a dense extension of a completely $0$-simple semigroup. 
In particular, it turns out that the structure theory of these rings is of a multiplicative nature since each isomorphism between the multiplicative semigroups of these rings is necessarily an isomorphism between the rings.
Let us also mention that the theory of dense extensions of semigroups was explored mainly by Russian authors (Ljapin, Gluskin, Shevrin). 
Mario is one among the few semigroup theorists from the Western world who have worked in this area.

Besides the account on this pair of papers, it should be mentioned that Mario Petrich has published a number of substantial papers on ideal extensions and the translational hull of semigroups and rings.
The translational hulls of the completely simple components play a significant role
in Lallement's structure theorem for completely regular semigroups mentioned above. 
Having at hand his description \cite{Pet1968} of the translational hull of a regular Rees matrix semigroup,
Petrich clears up in \cite{Pet1974} exactly how these components interact.
This construction of every completely regular semigroup from a semilattice, a family of Rees matrix semigroups and families of functions among their ingredients is refined in two further papers \cite{Cli-Pet1977} and \cite{Pet1987}, the former written jointly with Clifford (see also \cite[Theorem VI.5.2]{Gril1995}).
In \cite{Pet1969} Mario provides a fundamental faithful representation, called the Petrich representation in \cite[Theorem VI.2.3]{Gril1995}, for every regular semigroup $S$ by means of translational hulls of completely $0$-simple and completely simple semigroups.
The completely $0$-simple and simple semigroups involved are the traces of the $\mathcal{D}$-classes of $S$.
Let us call attention also to \cite{Pet1985} where the main tool is the translational hull of rings and, among others, a new access to Everett sums (these correspond to the Schreier extensions of groups) is presented.
In this connection it is also worth mentioning the survey paper \cite{Pet1970}, a very readable and extensive account on the topic up to the early seventies.

Another line of development in Mario's research inspired by a celebrated theorem in ring theory concerns semigroups of quotients and orders in semigroups.
Goldie's Theorem characterizes rings having a semisimple Artinian ring of quotients.
Early attempts to generalize it for semigroups assume that the semigroup has an identity element, and use the group of units to form quotients, but this has not led to significant results. Noticing that the role of the identity element in a ring is often taken over by the set of all idempotents in structural investigations into semigroups, Fountain persuaded Petrich to consider group inverses, which led to their pioneering paper \cite{FouPet1986}: the inverse of an element is taken in a subgroup around an arbitrary idempotent in the semigroup of quotients. 
They provided a list of necessary and sufficient conditions for a semigroup to have a completely $0$-simple semigroup of quotients in this sense. 
This result can be considered as a close analogue of Goldie's Theorem. 
This paper started a long line of investigations, due largely to Victoria Gould, whereby elements satisfying an additional condition (called square-cancellable, a notion introduced by Fountain) are required to have a group inverse in the semigroup of quotients. 
This condition is necessary when considering semigroups of quotients in general, in completely $0$-simple semigroups it is automatically satisfied.  
Mario had further contributions to the topic in cases where this condition is automatically satisfied. 
From these papers, let us quote just one result \cite{FouPet1996}: a semigroup is an order in a normal band of groups if and only if it is a normal band of reversible cancellative semigroups.

It has been noticed by Wagner \cite{Wag1953} that, similarly to groups where each congruence is determined by its class containing the identity element, each congruence of an inverse semigroup is determined by its system of classes containing idempotents.
A characterization of such systems is due to Preston \cite{Pre1954}, see also \cite[Theorem 7.48]{CliPre1961-7}. 
In \cite{ReiSch1967}, Reilly and Scheiblich study the pairs of congruences defining the same partition on the set of idempotents, and in \cite{Sch1974}, Scheiblich establishes that every inverse semigroup congruence $\rho$ is uniqueIy determined by the union of $\rho$-classes containing idempotents, which forms an inverse subsemigroup and is called the {\em kernel} of $\rho$, and by the restriction of $\rho$ to the semilattice of idempotents, called the {\em trace} of $\rho$. 
Inspired by these ideas and results, Mario \cite{Pet1978} introduces the notion of a {\em congruence pair} consisting of a full inverse subsemigroup and a congruence on the semilattice of idempotents, and presents necessary and sufficient conditions for actually building a congruence from a congruence pair.
One of the axioms of a congruence pair is further simplified in his book \cite{Pet1984}. 
The monographs \cite{How1995} and \cite{LawM1998} published since then prefer the kernel-trace description of inverse semigroup congruences.

In a joint paper by Petrich and Reilly \cite{PetRei1982}, kernel and trace are studied from a new perspective, providing unexpected generalizations. 
Given a congruence $\rho$ on an inverse semigroup $S$, the authors denote by 
$\rho_{\operatorname{min}}$ and $\rho^{\operatorname{min}}$ the least congruence on $S$ having the same trace and the same kernel as $\rho$, respectively.
Then they alternate applications of the two operators to the universal congruence, thus defining the {\em min network} of $S$.
The uppermost part of this partially ordered set of congruences features well-studied congruences, namely the least congruence such that the respective quotient semigroup is a group, semilattice, Clifford semigroup, $E$-unitary and $E$-reflexive semigroup, respectively.
It turns out that the class of inverse semigroups associated to each member of the min-network in this way forms a quasivariety, and a basis of quasi-identities is provided for each class.

The analogy between congruences of groups and inverse semigroups leads naturally to the idea of generalizing Schreier's theory of group extensions.
Mario \cite{Pet1981} introduces a general framework for treating the analogue of a group extension within the class of inverse semigroups, called a normal extension, and presents a solution for the normal extension problem in several special cases. 
The general solution of the normal extension problem for inverse semigroups is due to Allouch \cite{All1979}.
Pastijn and Petrich \cite{PasPet1986} generalize the kernel-trace approach to congruences on regular semigroups.
In \cite{PasPet1985} they introduce a general framework to study extensions in the class of regular semigroups, present in a number of cases how former results fit in this 
framework, 
and analyse the connections between former results.
In order to make it clear how much congruences and extensions might be more complicated for a regular semigroup than for an inverse semigroup, let us mention that neither the kernel nor the base set of the trace of a congruence on a regular semigroup forms a subsemigroup.
Let us also mention that operators arising from kernels, traces and one-sided versions of traces play a significant role in the description of the lattice $\mathcal{L}(\mathcal{CR})$ of all varieties of completely regular semigroups.

Mario's interest in completely regular semigroups has been an important aspect of his work for over fifty years.
The special form of his structure theorem, mentioned above from \cite{Pet1974}, for bands appeared earlier in the lecture notes \cite[Chapter 8]{Pet1967}, and played an important role in the unravelling of the lattice of band varieties and their bases of identities. 
An inequivalent system of band identities in up to three variables was announced by 
Kimura \cite{Kim1958}, and partial proofs were provided by Yamada \cite{Yam1962}.
Using his structure theorem for bands, Mario provided complete proofs and showed that 
to determine bases of identitites for all band varieties it would be sufficient to consider identities that conformed to one of five descriptions based on certain invariants of words.  
Building on his observations, Fennemore gave a complete description of the lattice of varieties of bands and bases of identities in his doctoral thesis \cite{Fen1969} and in
\cite{Fen1971}.
Similar results were obtained independently, and by completely different methods 
by Birjukov \cite{Bir1970} and Gerhard \cite{Ger1970}, the latter from his dissertation at McMaster University, 1968.

This line of development is pursued in two subsequent papers \cite{GerPet1986} and \cite{GerPet1989} by Gerhard and Petrich jointly. 
In the first paper the authors introduce an analogue of reduced words for groups. 
They define a new invariant $b(w)$ of any word $w$ on a set $X$ inductively 
such that the mapping $b$ solves the word problem for the free band on $X$, and the set of all words of the form $b(w)$ endowed with the 
multiplication $v\cdot w = b(vw)$ forms a model for it.
This approach is extended later by Mario and the first author in \cite{PetSil2000} and \cite{PetSil2002} for every relatively free band.
The second paper is something of a landmark.  
It introduces a number of new invariants for words, and uses them to accomplish both a new derivation of the lattice of band varieties and also the solution of the word problem for every such variety.  
The treatment is entirely semantic and liberated the development from the structure theorem.  
Three recursively defined systems of words 
$G_n,\,H_n,\,I_n\ (n\ge 2)$ 
are introduced, and it is shown that every proper subvariety of the variety $\mathcal{B}$ of all bands containing semilattices has a basis that consists of either one or two identities of the form 
\begin{equation} \label{GHI}
G_n = H_n,\ G_n = I_n,\ \overline{G_n} = \overline{H_n},\ \overline{G_n} = \overline{I_n}
\end{equation}
where $\overline{w} = x_k\cdots x_1$ for any word $w = x_1\cdots x_k$.  
In addition, those defined by  precisely one of these identities are 
exactly the join irreducible varieties of bands.  
This formulation of the identities was extremely fortuitous as these identities could be adapted to describe another very important family (see below).  

There are many classes $\mathcal{P}$ of completely regular semigroups for which the mapping 
$\theta_{\mathcal{P}}\colon \mathcal{V} \mapsto \mathcal{V}\cap\mathcal{P}\ 
(\mathcal{V} \in \mathcal{L}(\mathcal{CR}))$ 
induces an interesting equivalence relation on the lattice  $\mathcal{L}(\mathcal{CR})$ of all varieties of completely regular semigroups and
for which the classes are intervals
$[\mathcal{V}_{\theta_{\mathcal{P}}}, \mathcal{V}^{\theta_{\mathcal{P}}}]$.  
Such situations lead naturally to the operators 
$\mathcal{V} \mapsto \mathcal{V}_{\theta_{\mathcal{P}}}$ and 
$\mathcal{V} \mapsto \mathcal{V}^{\theta_{\mathcal{P}}}$ on $\mathcal{L}(\mathcal{CR})$.  
The investigation of these operators singly or in combinations 
reveals a great deal of information regarding the structure of $\mathcal{L}(\mathcal{CR})$, some of it local in nature and some of it more global.   
To pick an example, the operators induced by $\theta_{\mathcal{M}}$ are studied, among others, by Petrich and Reilly in \cite{PetRei1990}.
Here $\mathcal{M}$ denotes the class of monoids,
$\mathcal{V}_{\theta_{\mathcal{M}}}$ is the variety generated by 
$\mathcal{V}\cap\mathcal{M}$, and $\mathcal{V}^{\theta_{\mathcal{M}}}$, 
usually denoted $\mathcal{V}^L$, is the variety consisting of all completely regular semigroups all of whose submonoids belong to $\mathcal{V}$.
Another class of completely regular semigroups providing a significant operator in a similar way is $\mathcal{B}$. 
It is due to Trotter \cite{Tro1989} that the equivalence classes induced by $\theta_{\mathcal{B}}$ are intervals.
Clearly $\mathcal{V}_{\theta_{\mathcal{B}}}=\mathcal{V}\cap\mathcal{B}$, and
the variety $\mathcal{V}^{\theta_{\mathcal{B}}}$ is usually denoted  $\mathcal{V}^B$.

Mario has made a very significant contribution to literature on various operators on $\mathcal{L}(\mathcal{CR})$: 
kernel $K$, trace $T$, left trace $T_{\ell}$, left kernel $K_{\ell}$, core $C$, local $L$, bands $B$.  
He is a definite leader when it comes to analysing the properties of combinations of operators, especially $K, T, T_{\ell}$ with $C$ and $L$, see for example \cite{Pet2018}.
These excursions encounter a multitude of varieties for which, with very few exceptions,  he provides a basis of identities. 

In \cite{Pet2007}, Mario introduces a generalization of the systems of the words 
$G_n,\,H_n,\,I_n\ (n\ge 2)$ (see (\ref{GHI})) in the signature of completely regular semigroups, and calls the varieties of left and right normal orthogroups together with the varieties defined by any one of the generalized versions of the identities in (\ref{GHI})
{\em canonical varieties}.
He establishes for any band variety $\mathcal{V}$
containing the variety of semilattices that $\mathcal{V}^B$ is equal to the intersection of all the canonical varieties containing $\mathcal{V}$.  Consequently, 
$\mathcal{V} = \mathcal{V}^B$ for every canonical variety $\mathcal{V}$, and every variety 
containing the variety of semilattices and of the form $\mathcal{V}^B$ is an intersection of canonical varieties.    
In a series of interesting papers \cite{Pet2014}--\cite{Pet2017} and \cite{Pet2019}, he then investigates the sublattice $\Gamma$ of $\mathcal{L}(\mathcal{CR})$ generated by the canonical varieties.  
In each of these papers, he provides a precise picture of successively larger parts of $\Gamma$ together with a great deal of information regarding both the general structure of the lattice, including $T_{\ell}, T_r, C$ and $L$-classes, as well as specific properties of a good few of the varieties that appear low down in the lattice.  
Clearly, $\Gamma$ is infinite but, what is more suprising, is that it is distributive and of width at most four.  
At this point, there is no complete description of $\Gamma$.  However, in the last cited paper, Mario provides a remarkable conjecture regarding its structure. 
He provides sufficient information about $\Gamma$ to establish certain features of the conjecture, such as an unexpected layered structure with repeated patterns.  
There is no reason to believe that he has, as yet, laid down his quill and perhaps a confirmation or otherwise of this conjecture could still be in preparation.

Besides his first book {\em Introduction to Semigroups} \cite{Pet1973m} and the research notes {\em Regular Semigroups as Extensions} \cite{PasPet1985} mentioned above in some detail, Mario has published six monographs, one of them with a co-author.
All but two focus on a field of semigroup theory, the rest being about interrelations between {\em Rings and Semigroups} \cite{Pet1974m} and about {\em Categories of Algebraic Systems} \cite{Pet1976m}.
His most cited monographs are {\em Inverse Semigroups} \cite{Pet1984} and {\em Completely Regular Semigroups} \cite{PetRei1999}, the latter co-authored by Reilly.
The first one, published in 1984, was the first reference book on the subject, and continues to be abundantly cited by researchers in the area. 
It contains most of the significant results known at that time and, simultaneously, it is a self-contained and carefully organized text that has helped young researchers to deepen their knowledge on the subject.
The second one also belongs on the shelf of anyone interested in regular semigroups.
Similarly to the first one, it is a comprehensive introduction to the theory of semigroups in the title with exercises at the ends of the chapters for newcomers.
It is sad that the most highlighted results of the subject --- the descriptions of the free completely regular semigroups and the lattice of all varieties of completely regular semigroups --- are postponed to another volume.
Many researchers still hope that this volume will be published. 
A much earlier monograph {\em Structure of Regular Semigroups} \cite{:Pet1977r} contains results on classes of regular semigroups that were the focus of his attention in the 1970's.
Finally, the one entitled {\em Lectures in semigroups} \cite{:Pet1977l} focuses on bands, on decompositions of semigroups into normal bands of semigroups and on semigroups whose subsemigroup lattices are modular, distributive or complemented.

\medskip

\begin{center}
	\bf Acknowledgement
\end{center}

The authors express their gratitude to Norman Reilly for his valuable advice and generous help during the preparation of this tribute.
Their thanks go also to L\'{a}szl\'{o} M\'{a}rki for his useful suggestions.

The first author was partially supported by CMUP, which is financed by national funds through FCT -- Funda\c c\~ao para a Ci\^encia e a Tecnologia, I.P., under the project with reference UIDB/00144/2020,
and the second author by the National Research, Development and Innovation Office (Hungary), grant no.\ K128042.

\end{document}